\documentclass[12pt]{article}
\usepackage{amsmath,amssymb,latexsym,graphicx}
\newtheorem{theorem}{Theorem}
\newtheorem{lemma}[theorem]{Lemma}
\newtheorem{proposition}[theorem]{Proposition}
\newtheorem{remark}[theorem]{Remark}
\newtheorem{corollary}[theorem]{Corollary}
\newtheorem{definition}[theorem]{Definition}

\def\pf{{\bf Proof }}

\begin{document}
\title{Characterization of polyanalytic functions by meromorphic extensions into chains of  circles} 
\author{Mark~ L.~Agranovsky}
\maketitle

\begin{abstract}
Let $C_t=\{ z \in \mathbb C:|z-c(t)|=r(t), t \in (0,1)\}$ be a $C^1$-family of circles in the plane, such that $\lim_{t \to 0+} C_t=\{a\}, \ \lim_{t \to 1-}C_t=\{b\}, \ a \neq b.$
The discriminant set $S$ of the family is defined as the closure of the set 
$\{c(t)+r(t)w(t), t \in [0,1]\}$, where $w=w(t)$ is the root of the quadratic equation 
$\overline c^{\prime}(t)w^2+2r^{\prime}(t)w+c^{\prime}(t)=0$  with $|w|<1$ (if such a root exists) and
$|c^{\prime}(t)|^2+|r^{\prime}(t)|^2 \neq 0.$
Suppose that (*) $S$ does not contain a continuous curve joining $a$ and $b$. We prove that if $f \in C^{\nu}(\Omega),$ \ $\Omega=\cup C_t,$ is a regular, in a certain sense, function and $f$ possesses, for each $t \in (0,1),$ a meromorphic extension inside $C_t$ with the only singular point-a pole at $c(t)$ of order at most $\nu \in \{0\} \cup N,$
then $f$ is polyanalytic of order $\nu$, i.e., $f(z)=h_0(z)+h_1(z)\overline z+\cdots \overline z^{\nu}h_{\nu}(z),$ 
where the $h_j's$ are analytic functions in $\Omega.$    
For $\nu=0$ the condition (*) can be omitted. If $|c^{\prime}(t)| >|r^{\prime}(t)|, t \in (0,1),$ then $S=\emptyset$ and
the condition of regularity can be dropped. A hyperbolic version of the result is given.
\end{abstract}

\section{Introduction}\label{S:Intro}
The problem of testing analyticity on families of closed curves was studied in many articles.
Generally speaking, the problem can be described as follows: suppose that a domain $\Omega \subset \mathbb C$ is covered
by a family of closed curves, must a continuous function $f$ be analytic in $\Omega$ if it extends analytically inside each curve?

In the paper \cite{AV} the case of translation and rotation invariant families of curves in the plane was solved. Some special cases  were considered in \cite{GG}, \cite{E}.
In \cite{G1, G2, G3}
the case of rotation invariant families of curves was investigated in details. In \cite{AG}
the case of arbitrary continuous families of circles was solved for rational (in the real variables) functions.
The case of continuous functions and pretty large families of circles was solved in \cite{T1,T2}. In general setting, for generic families of Jordan curves, the problem was solved (in real-analytic category)  in \cite{A0,A1}.

This article is devoted to similar characterization of classes of functions, larger than analytic, namely, of polyanalytic functions. These functions are characterized by the condition of meromorphic extendibility
rather than by analytic that.

Our testing families are pretty general smooth families of circles in the plane, arising from one point and shrinking to another one.  We will call such families {\it chains.}  In fact, the condition of shrinking to two points can be replaced by periodicity of the family, or, presumably, for non-periodic families, by the condition of existing  two circles belonging to the exterior of each other.
In the recent article \cite{G4} test of analyticity in the unit disc was proved for a special chain of the above type: union of  circles of radii $\leq 1$ with the centers at the origin,  and horicycles through a fixed point (circles tangent from inside to the unit circle at a fixed point).

The main results is the following. We describe a large family of chains of circles such that
the  meromorphic extendibility into each circle from this chain, with the pole at the center,
implies that the function is polyanalytic. Of course, the result includes, as a particular case,
characterization of analytic functions.

The interest to the problem considered is this article is, in particular, related to an application to a problem in several complex variables.  These applications will be given elsewhere.   

\section{Main results}\label{S:main}

We will use the notations 
$\Delta=\{z \in \mathbb C: |z|<1\}, \ C(c,r)=\{z \in \mathbb C:|z-c|=r\}, \ D(c,r)=\{z \in \mathbb C:|z-c| <r\}.$

\subsection{Formulation of the main results}
Let us formulate precisely the main result of the article. First of all, some  definitions.
Let 
$$c:[0,1] \to \mathbb C,  \  r:[0,1] \to [0,\infty)$$ 
be two mappings belonging to the class $C[0,1] \cap C^1(0,1).$
We assume that 
$$|c^{\prime}(t)|^2+|r^{\prime}(t)|^2 \neq 0,$$
except, maybe, for  a  finite set $t_1, \cdots, t_k \in (0,1),$ and the zero at each point $t_j$ is of a finite order.
 
We will call the family $\mathcal C=\{C_t\}_{t \in [0,1]}$ of the circles
$$C_t=C(a(t),r(t)):=\{z \in \mathbb C: |z-c(t)|=r(t)\}$$ a {\it chain of circles}, with the initial point $a$ and the end point $b,$
if
$$C_0=\{a\} \ \mbox{and} \ C_1=\{b\}, \ \mbox{i.e.}, \ r(0)=r(1)=0, \ c(0)=a, \ c(1)=b.$$
We will denote also $D_t=D(c(t),r(t)).$

Associate  with the chain $\mathcal C=\{C_t, \ t \in [0,1]\}$  a set constructed in the following way. 
Consider the quadratic polynomial
\begin{equation}\label{E:quadratic}
d(w,t):=\overline c^{\prime}(t)w^2+2r^{\prime}(t)w+c^{\prime}(t)
\end{equation}
We will call $d(w,t)$ the {\it discriminant} of the chain $\mathcal C.$ 

Suppose that $t$ is such that $d(w,t)$ is not identical zero and let $w_1,w_2$ be the roots of the discriminant $d(w,t).$ If $c^{\prime}=0$ then $w_1=w_2=0.$ Otherwise,  $|w_1||w_2|=1,$  and  either both roots lie on the unit circle, which happens if and only if $|c^{\prime}| \geq |r^{\prime}(t)|$, or one of them is inside the unit circle and another one is outside. 
Define the set $S=S(\mathcal C)$ as $$S=\mbox{closure}\{c(t)+r(t)w: d(w,t)=0, |c^{\prime}(t)|^2+r^{\prime}(t)|^2| \neq 0, \ |w|<1, \ t \in [0,1]\}.$$
We will call $S(\mathcal C)$  the {\it discriminant set} associated with the chain $\mathcal C.$

Let $\Omega$ be a closed domain in $\mathbb C$ with smooth boundary. We say that a function $f \in C(\Omega)$ 
is {\it regular} if either $f=0$ or the zero set of $f$ consists of curves and points and $f$ vanishes on them to finite order . The definition in more precise form is given in Section \ref{S:traveling}.  
An example of such function are delivered by real-analytic functions in $\Omega.$

Let $\Gamma=\{\gamma_t\}, \ t \in (0,1)$ be a one-parameter family of differentiable Jordan curves in the plane. Denote $D_t$ the domain bounded by the curve $\gamma_t.$ 

\begin{definition}\label{D:merom}
Let $g$ be a function defined on the union
$\Omega$ of the curves $\gamma_t.$ 
We say that $g$ meromorphically extends from  $\Gamma$ if for every $t \in (0,1)$ there exists a function $G_t,$ meromorphic in the domain $D_t$ and continuous in the closure $\overline D_t$, such that $G_t(z)=g (z)$ for $z \in \gamma_t.$
\end{definition}

The main result of this article is the following

\begin{theorem}\label{T:MAIN}
Let $\mathcal C=\{C_t=C(c(t),r(t)), \ t \in [0,1]\}$  
be a chain of circles with the initial and end points $a$ and $b$, ($a \neq b$), 
correspondingly. Suppose  that 

(*) the points $a$ and $b$ cannot be joined by a continuous curve within
the discriminant set $S(\mathcal C).$ 

Let $\nu \in \{0\} \cup \mathbb N.$ If $f \in C^{\nu}(\Omega)$ is a  regular function in $\Omega:=\cup_{t \in [0,1]}C_t,$  
such that for any $t \in (0,1)$
the restriction $f\vert_{C_t}$ extends in the disc $D_t$ as a meromorphic function with the only singularity - a pole of the order at most $\nu$ at the center $c(t)$, then $f$ is polyanalytic function of order $\nu$: 
$f(z)=h_0(z)+\overline z h_1(z) +\cdots + \overline z^{\nu}h_{\nu}(z).$ 
In this case, $f$ possesses the above type of meromorphic extension inside any circle in $\Omega.$
In the case $\nu=0$ the condition (*) can be omitted.
\end{theorem}
When $\nu=0$ then Theorem \ref{T:MAIN} reduces to a  test for analytic functions.
Polyanalytic functions of order $\nu$ are solutions of the  the equation 
$$\frac{\partial^{\nu+1}f}{\partial \overline{z}^{\nu+1}} =0.$$
Theorem \ref{T:MAIN} characterizes solutions of this PDE  in terms of meromorphic extensions into one-parameter chains of circles. 

Stronger conditions for the chain of circles allow to omit regularity conditions for functions: 
\begin{theorem}\label{T:noregular} Suppose that in Theorem \ref{T:MAIN} the additional conditions hold: 
$|c^{\prime}(t)|> |r^{\prime}(t)|, \ t \in (0,1)$  and none of two circles $C_t, \ C_s, \ t \neq s,$ is contained strictly inside one another. Then the condition of regularity of $f$ in Theorem \ref{T:MAIN} can be dropped. 
\end{theorem} 
The condition $|c^{\prime}(t)|>|r^{\prime}(t)|$ first appeared in \cite{AG} and then in \cite{T2}.

Notice, that the families of circles, involved in Theorem \ref{T:MAIN}, are of more general type then 
those in  \cite{T2}, where the result is proved for families of the type considered in Theorem \ref{T:noregular}. 
Hence, even for $\nu=0$ (characterization of analytic functions) the result formulated in Theorem \ref{T:MAIN} is new. 
Moreover, our main examples correspond just to the case of enclosed circles, because this case is important in applications.    

\subsection{Examples and corollaries}\label{S:examples}
Let us give some concrete examples of chains for which Theorem \ref{T:MAIN} works. None of the chains
in these examples belongs to the type covered by Theorem \ref{T:noregular}.

\noindent
1.{\it Hyperbolic circles}.

For each $a \in \Delta$, denote
$$H(a,r)=\{z \in \mathbb C: |\frac{z-a}{1-\overline a z}|=r\}$$
the hyperbolic circle in the unit disc $\Delta$, centered at $a$.  
The geometric center of the circle $H(a,r)$ is
\begin{equation}\label{E:geometric}
c=c(a,r)=a\frac{1-r^2}{1-|a|^2r^2}.
\end{equation}
When $a=0$ then the hyperbolic center coincides with the geometric one, and also $H(a,1)=\partial \Delta$ is the unit circle.

Let  a function $r \in C[0,1] \cap C^{1}(0,1)$ be such that $0 \leq r(t) \leq 1, \ t \in [0,1];
\ r(t) > 0, \ t \in (0,1); \ r(0)=r(1)=0, \ r(1/2)=1.$ We also assume that $r^{\prime}(t)=0$ 
only at finite number of points. For instance, we can choose $r(t)$ strictly  increasing on $[0,1/2)$ and strictly decreasing on $(1/2,1]$ so that $r^{\prime}(t)=0$  only at $t=1/2.$  
An example of such function is $r(t)=4t(1-t).$

Given two distinct points $a,b \in \ \Delta$, define the chain of circles $\mathcal H(a,b)=\{C_t\}_{t \in [0,1]}$ by:
\begin{equation} \label{E:hyper}
C_t=H(a,r(t)), \ t \in [0,\frac{1}{2}]; \ \ C_t=H(b,r(t)), \ t \in [\frac{1}{2},1].
\end{equation}
The circles in $\mathcal C=\{C_t\}_{t \in [0,1]}$ arise from the point $a$ and grow, keeping the  hyperbolic center $a,$ till they reach the unit circle at the moment $t=1/2.$ After then the circles switch the hyperbolic centers from $a$ to $b$ and start shrinking till they degenerate to the point $b$.
It is checked in Section \ref{S:proof2} that $\mathcal C$  satisfies the condition (*) from Theorem \ref{T:MAIN}.

\noindent
2. {\it Horicycles.}

Let $a \in \partial \Delta.$ Horicycles are the circles $Hor(a,r)$ of radius $r \leq 1,$ containing $a$ and entirely
belonging to  $\overline \Delta.$ Clearly, $O(a,1)=\partial \Delta.$
For any radius-function $r(t)$ from the previous example, define
\begin{equation} \label{E:hori}
C_t=Hor(a,r(t)), \ t \in [0,\frac{1}{2}]; \ \ C_t=Hor(b,r(t)), \ t \in [\frac{1}{2},1].
\end{equation}
This is a smooth chain of circles covering the closed domain $\Omega=\overline \Delta$. 
The condition (*) can be checked immediately in this case.
Indeed, take for simplicity $a=-1, b =1$ Consider $t$ in the half-interval $[0,1/2).$ 
Then the center $c(t)=-1+r(t).$  The equation (\ref{E:quadratic})
reads as $r^{\prime}w^2+2r^{\prime}w+r^{\prime}=0$ and since $r^{\prime}(t) \neq 0$ then there is one unmovable
double root $w(t)=-1.$ Then $c(t)+r(t)w(t)=(-1+r)+r(-1)=-1.$ Analogously $c(t)+r(t)w(t)=1, \ t \in (1/2,1].$
Thus the discriminant set $S(\mathcal C)$ consists of two points $a=-1$ and $b=1$ and the condition (*) is fullfiled for the obvious reason. Therefore, this is the case covered by Theorem \ref{T:MAIN}.

\noindent
3. {\it A chain of a mixed type: concentric circles and horicycles}.

The following family of a mixed type was considered in \cite{G4}.
\begin{equation} \label{E:hyper_hori}
C_t=C(0,r(t)), \ t \in [0,\frac{1}{2}]; \ \ C_t=Hor(b,r(t)), \ t \in [\frac{1}{2},1].
\end{equation}
Here $|a|<1$ and $|b|=1.$
The first part of this family consists of circles centered at the origin and the second one-of horicycles
through the boundary point $b \in \partial \Delta$. The two families are joined through the unit circle, corresponding  to the value $1/2$ of the parameter $t$. This a chain of circles, with the initial point $a=0$ and the end point $b$ and covering the domain $\Omega=\overline \Delta$. 

In \cite{G4} the analyticity test for such family of circles was proved. 
It is easy to check that the condition (*) holds for this chain. 
Indeed, as we have seen in the previous example, for the horicycle part of the family we have
$c(t)+r(t)w(t)=b.$ On the other hand, for the first part, 
consisting of Euclidean circles with the centers at the origin, we have $c(t)=0,$ the equation (\ref{E:quadratic}) has the solution $w=0$ and hence $c(t)+r(t)w=0$. Thus, the discriminant set $S(\mathcal C)$ again consists of two points $a=0$ and $b$ and these point can not be joined by a continuous curve in $S(\mathcal C).$
Therefore, the analyticity test of \cite{G4} is a particular case $\nu=0$ of Theorem \ref{T:MAIN} for a special family of circles. However, the result of \cite{G4} is true for continuous functions, while we impose the stronger condition of regularity.

We will finish this section with a  version of Theorem \ref{T:MAIN} for the special case considered in Example 1. 
In fact, this case was a starting point and a motivation for this article, due to its applications which we will present elsewhere.

Let $F$ be a function in $\Delta$. We will write
$$F =O(\frac{1}{(1-|z|^2)^{\nu}}), \ |z| \to 1,$$
if $(1-|z|^2)^{\nu}F(z)$ continuously extends to the closed disc $\overline \Delta$ and the extension has on the
unit circle $\partial \Delta$ only isolated zeros of finite order.

\begin{corollary} \label{C:hyperbolic}
Let $F \in C^{\nu}(\overline \Delta)$ be a regular function and
$F(z)=O((1-|z|^2)^{-\nu}), \ |z| \to 1$. Suppose that for some $a,b \in \Delta, \ a \neq b,$ the function $f$
extends meromorphically from the hyperbolic circles $H(a,r), H(b,r), \ r \in (0,1],$ with the only singular point-
a pole of the order at most $\nu$ at the hyperbolic centers $a$ and $b,$ correspondingly.
Then $F$ has the form
\begin{equation}\label{E:F}
F(z)=h_0(z)+\frac{h_1(z)}{1-|z|^2}+ \cdots +\frac{h_{\nu}(z)}{(1-|z|^2)^{\nu}},
\end{equation}
where $h_j(z), \ j=0,\cdots,\nu$ are analytic functions in $\Delta.$
\end {corollary} 
  
\begin{remark}\label{integral}
The condition of meromorphic extendibility inside the circles can be 
written in integral form, as vanishing certain 
complex moments. In Theorem \ref{T:MAIN} and Theorem \ref{T:noregular} the  condition reads as follows:
$$\int\limits_{|z-c(t)|=r(t)}f(z)(z-c(t))^m dz=0 , t \in [0,1], \  m \geq \nu.$$
In Corollary \ref{C:hyperbolic} it can be written as :
$$\int\limits_{|z|=r}f(\frac{z+a}{1+\overline a z})z^m dz=0, \ r \in [0,1], \  
\int\limits_{|z|=r} f(\frac{z+b}{1+\overline b z})z^m dw=0, \ r \in [0,1], \ m \geq \nu.$$
Therefore, the above formulated results can be regarded as theorems of Morera type.
\end{remark}
\subsection{Comments on the proofs and plan of the article} 
The proof of Theorem \ref{T:MAIN} is essentially  based on ideas and methods from \cite{A0,A1,A2}, developing argument principle for parametric families of holomorphic mappings. 

The main idea is to study the dynamics of zeros and poles of the meromorphic extensions
of $f$ into the circles $C_t$. 
Due to the regularity condition for $f$ this dynamics is nice and can be well understood. 
When the circles $C_t$ move from $a$ to $b,$ the zeros and poles move along certain curves.
We introduce the notion of traveling zeros and poles for those which fill a continuous path from $a$ to $b$.

The key tool is Proposition \ref{P:zeros_poles} proved in Section \ref{S:traveling}. 
It concerns general one-parameter smooth chains of closed curves, not only circles, and states that for nonzero regular functions the number of traveling zeros equals to the number of traveling poles. 

Then the proof of Theorem \ref{T:MAIN} goes by the induction in the order $\nu$ of poles of the meromorphic extensions at the centers $c(t)$ of the circles $C_t$, under assumption that there is no other traveling poles.

The beginning of the induction and the induction step are done by applying Proposition \ref{P:zeros_poles}
to the $\overline \partial$-derivatives. Namely, we prove  in Section \ref{S:dbar} 
that the differentiation in $\overline z$ reduces the order $\nu>0$ of the poles at the centers $c(t),$ 
albeit new simple pole inside $C_t$  may appear at the discriminant set $S(\mathcal C)$.
However, the condition (*) guarantees that this discriminant set does not contain continuous curves from $a$ to $b.$ 
This means $S(\mathcal C)$ do not contribute to the number of traveling poles and hence no new traveling poles different from those at $c(t)$ appear. Then by the  assumption of the induction we have that $\partial_{\overline z} f$ is polyanalytic of order $\nu-1$ and this is equivalent to $f$ being polyanalytic of order $\nu.$ 

The beginning of the induction  corresponds to $\nu=0$, i.e. to the case when the function  extends analytically inside the circles $C_t.$ Then we show that the function
 $g:=\overline\partial_{\overline z}f$ does the same and, moreover, its meromorphic extension has zero at $c(t)$ of at least second order. Therefore, even if the discriminant set $S(\mathcal C)$ contributes a simple traveling pole, still, the
number of traveling zeros is bigger. Then Proposition \ref{P:zeros_poles} implies 
$\overline\partial_{\overline z}f=0$, i.e. $f$ is analytic. Notice that the condition (*) for the discriminant set is not needed in this case. The final part of the
proof of Theorem \ref{T:MAIN} is presented in Section \ref{S:proof1}.

Theorem \ref{T:noregular} is easier. In this case, the induction step goes the same line, but this time there is no 
need in applying Proposition $\ref{P:zeros_poles}$
because the discriminant set is empty and hence $\overline z$- differentiation does not produce  poles out of the centers $c(t).$ As for the initial step of the induction is concerned, it is provided  by the result of Tumanov \cite{T2} which is true just for chains like in Theorem \ref{T:noregular} and  requires only continuity of the functions.

In Section \ref{S:proof2} we prove Corollary \ref{C:hyperbolic}. Section \ref {S:concluding} is devoted to concluding remarks.
\section{Dynamics of zeros and poles}\label{S:traveling}
In this section we prove a key statement on meromorphic extensions into parametric families of Jordan curves. This statement concerns general curves, not necessarily circles,  and has
meaning of argument principle for parametric families of closed curves. The proof is essentially  based on methods  of  \cite{A1,A2}.

We start with some definitions.

Let $\Gamma=\{\gamma_t\}, \ t \in (0,1),$ be a one-parameter family of differentiable Jordan curves in the plane. Denote $D_t$ the domain bounded by the curve $\gamma_t.$ 

Suppose that the curves $\gamma_t$ depend continuously on the parameter $t \in (0,1)$. Suppose that a function $g$ meromorphically extends from $\Gamma =\{\gamma_t\}$ (see Definition \ref{D:merom}) and $G_t$ are the corresponding meromorphic extensions . 

\begin{definition}\label{D:traveling}
We say that a continuous curve
$Z  \subset \Omega$ is a {\it traveling zero} of the meromorphic extensions of $g$ from the family $\Gamma$ if there exists a parametrization $z=z(t), \ t \in (0,1)$ such that $z(t) \in \overline D_t$ and $G_t(z(t))=0.$

We say that a continuous curve $P \subset \Omega$ is a {\it traveling pole} for $g$ if there exists a  parametrization 
$z=p(t) $ of the curve $P$ such that $p(t) \in  D_t$ is a pole of the function $G_t, \ t \in (0,1).$
\end{definition}

Let $\gamma$ be a differentiable curve in the plane and let $\rho$ be a defining function, i.e.
$\gamma=\{\rho=0\}$ and $\nabla \rho \neq 0$ on $\gamma.$ We say that a function $g$ in a neighborhood of $\gamma$
{\it vanishes on $\gamma$ to finite order } $m$ if $g=\rho^m g_0$ near $\gamma$, where $g_0$ has on $\gamma$ only isolated zeros of finite order.     
\begin{definition}\label{D:regular}
Let $\Omega$ be a bounded domain in $\mathbb C$ with smooth boundary. We say that the function $g \in C( \Omega)$ is regular if the zero set $g^{-1}(0)$ of $g$ consists of  isolated zeros of finite order of smooth curves in $\Omega$ on each of which $g$ vanishes to a finite order. 
\end{definition}

Examples of regular functions are nonzero real-analytic function in $ \Omega$.

Now we will prove our key Proposition. Before formulating it, let us consider an  example which illustrates
the statement we are going to prove.

\bigskip
\noindent
{\bf Example}
Consider a chain $C_t=C(c(t),r(t)), \ t \in [0,1],$ of circles shrinking to two distinct points:
$a=C_0, \ b=C_1.$ Let $f(z)=\overline z.$ This function has the following meromorphic extension into
any circle $C_t$: 
$$F_t(z)=\overline c(t)+ \frac{r^2(t)}{z-c(t)}.$$
This extension has the simple traveling pole $p(t)$ at the center $p(t)=c(t).$
The function $F_t$ has the simple zero at the point 
$$z_t=c(t)(1-\frac{r^2(t)}{|c(t)|^2}).$$
This zero is strictly inside the circle $C_t$ so long as $|c(t)| \geq r(t).$ When $|c(t)|=r(t)$
then $z_t=0.$ If $|c(t)|<r(t)$ then $z_t$ is outside of the circle $C_t.$
It follows that the continuous curve given by the parametrization
$$z(t)=c(t)(1-\frac{r^2(t)}{|c(t)|^2})_+$$
is a traveling zero for meromorphic extensions of $f(z).$
Thus, we see that the function $f(z)=\overline z$ has one traveling pole and one traveling zero. 

The following Proposition, which is a version of the argument principle for parametric families of domains (see \cite{A3}), says that this is always the case: the meromorphic extensions in chains have as many traveling poles as traveling zeros:

\begin{proposition}\label{P:zeros_poles}
Let $\Gamma=\{\gamma_t,\} \ t \in (0,1)$ be a $C^1$-family of differentiable Jordan curves, arising from a point $a \in \Delta$ ans shrinking t a point $b \in \Delta, a \neq b$, i.e.
$\lim\limits_{t \to 0+}  \gamma_t=\{a\}, \  \lim \limits_{t \to 1-}=\{b\}.$
Let $g \in C^1(\overline \Delta)$ be  a  regular function. Suppose that $g$ meromorphically extends from $\Gamma$ and
the meromorphic extensions have $N=N_g$ traveling zeros and $M=M_g$ traveling poles. Then $N_g=M_g.$

\end{proposition}
\pf  The proof follows the methods from \cite{A1, A2}. 

We start with the parameterizing the family $D_t$ of the domains, bounded by the curves $\gamma_t,$
by a smooth family of conformal mappings:
$$\omega_t: \Delta \mapsto D_t.$$
Define
$$\omega(\zeta,t):=\omega_t(\zeta).$$
Let $G_t$ be the meromorphic extension of $g\vert_{\gamma_t}$ into $D_t.$ Define
$$\varphi(\zeta,t): =G_t(\omega_t(\zeta)).$$
Thus, we have constructed two functions, $\omega$ and $\Phi$, defined in the solid cylinder $\overline \Delta \times [0,1].$ 
The function $\omega(\zeta,t)$ is analytic in $\zeta \in \Delta$ and the function $G(\zeta,t)$ is meromorphic in $\zeta$.

The zero sets $Z(\omega-w)=\{\omega(\zeta,t)-w=0\}, \ Z(\varphi)=\{\varphi(\zeta,t)=0\},$ as well as the set of the poles $P(\varphi)=\{\varphi(\zeta,t)=\infty\}$ of the function $\varphi$, consist of finite number of differentiable curves.

The proof of Proposition \ref{P:zeros_poles} is based on computing the following integral:

\begin{equation}\label{E:defI}
I(q)=\frac{1}{2\pi i} \int\limits_{\partial \Delta \times (0,1)}\frac{d\varphi}{\varphi} \wedge \frac {d\omega}{\omega -q},
\end{equation}

We will show that, on one hand, this integral is well defined and equals zero for 
$q \notin \overline \Omega$ and, on the other hand, it can be expressed in terms of
zeros and poles of the function $\varphi$.  Comparing these two expressions for $I(q)$ will lead to the equality between the numbers of traveling zeros and traveling poles of $\varphi.$. 

We start wit the first statement.
\begin{lemma}\label{L:I(q)=0}
Let $q \in \mathbb C.$  
The integral $I(q)$ is well defined for $q \notin \overline \Omega$ and 
moreover $I(q)=0$ for such $q.$
\end{lemma}
\pf

The change of the variable of integration $ z=\omega(\zeta,t)$ yields
\begin{equation}\label{E:change}
I(q)=\deg \varphi \int\limits_{\Omega} \frac{dg}{g} \wedge \frac{dz}{z-q}=deg \varphi
\int\limits_{\Omega} \frac{\frac{\partial g}{\partial \overline z}}{g}\frac{1}{z-q}d\overline z \wedge dz.
\end{equation}
where $\deg \varphi$ is the  topological (Brouwer) degree of the mapping 
$$\varphi: (\partial \Delta  \times (0,1)) \cup  (\overline \Delta \times (\{0\} \cup \{1\})) \mapsto 
\overline \Omega.$$ 
To prove that $I(q)$ is well defined, it suffices to prove  that the singular integral in the right hand side of (\ref{E:change}) converges. 
Since $ q \notin \overline \Omega$, all the singularities come from the zeros of $g$.
Due to the compactness of $\overline \Omega$, it suffices to prove only local integrability.

If $z_0$ is an isolated zero of $g$ then $g(z)=O(|z-z_0|^k),$ when $z \to z_0$ and then
$$\frac{g_{\overline z}}{g} =O(\frac{1}{|z-z_0|}), z \to z_0,$$ 
which in dimension 2 is an integrable singularity.

If $z_0$ is a non-isolated zero then, first,  by the condition, it can not be a boundary point 
of $\Omega$ and, second, $g$ vanishes on a smooth curve $L$ near $z_0.$ After applying a suitable local diffeomorphism, the curve $L$ 
can be written in local coordinates $(u_1,u_2)$ near $z_0$ as  $L=\{u_1=0\}.$ Then we have $g(u_1,u_2)=u_1^m g_0(u_1,u_2),$ near $z_0,$ for some natural $m$ and some function $g_0$ having only isolated zeros. Therefore,
$$\frac{g^{\prime}_{\overline z}}{g}=m\frac{1}{u_1}g_0+ \frac{(g_0)^{\prime}_{\overline z}}{g_0}.$$
The first term is integrable against $u_1,u_2$ in the sense of principal value, the second one 
has isolated singularities and is integrable as well.

Thus, the right hand side in (\ref{E:change}) is well-defined and to prove that $I(q)=0$ it suffices
to understand that $\deg \varphi=0.$  This is the case, because  $\varphi$ maps the closed manifold 
(the cylinder $\partial \Delta \times (0,1)$ completed by the "`bottom"' and "`top"' discs 
$\overline \Delta_ -=\overline \Delta \times \{0\}$ and $\overline\Delta_+= \overline \Delta \times \{1\}$) onto 
the manifold $\Omega$  which is a domain in $\mathbb C$ and  has the nonempty boundary.  

Albeit the domain of definition of $\varphi$  is not everywhere smooth, as it
has two edges, namely, the circles $\partial \overline\Delta_{\pm},$  one can approximate this manifold by smooth ones by smoothening  near the edges and this slight smoothening does not affect on the topological degree.This remark completes the proof.
\begin{lemma}\label{L:I(q)=}
Let $q \notin \overline \Omega$.  The integral $I(q)$ equals
\begin{equation}\label{E:I(q)=}
I(q)=\int\limits_{\omega(Z(\varphi))}\frac{dz}{z-q}-\int\limits_{\omega (P(\varphi))}\frac{dz}{z-q}.
\end{equation}
Here $\omega(Z(\varphi))$ and $\omega(P(\varphi))$ are understood as 1-chains, i.e. counting 
multiplicities of zeros or poles of $\varphi (\zeta,t),$  correspondingly, equipped by the orientation  
defined by the increasing parameter $t$.
\end{lemma}
\pf
Identity (\ref{E:I(q)=}) follows from the Stokes formula via technology of computing currents: 
we delete tubular neighborhoods of the zeros sets $Z(\varphi) \cup P(\varphi)$ and then, using
that the differential form in the integral (\ref{E:I(q)=}) is closed out of the singular set, replace
the integral by a  surface integral. The last step is shrinking the tubular neighborhoods to the sets
$Z(\varphi)$ and $P(\varphi)$ (see \cite{A1} for details).

However, there is an alternative  way of deriving the expression (\ref{E:I(q)=}), more in spirit of function theory. Namely, the identity (\ref{E:I(q)=}) can be obtained from the logarithmic residue formula applying in the variable $\zeta.$   
Of course, the proof of latter formula also is based on the Stokes formula, so that the alternative way  corresponds just to writing $I(q)$ as a double integral, in $zeta$ and $t$, applying Stokes formula in $\zeta$ and then integrating in $t.$ 

Let us give this computation in more detailed way.
The natural coordinates  of the points $ (e^{\psi},t) \in \partial \Delta \times (0,1)$ are the angle
$\psi$ and the parameter $t$.

The differential form $\eta$ under the sign of integral in  \ref{E:defI} writes  in these coordinates as

$$\eta:= \frac{1}{2\pi i}\frac{d\varphi}{\varphi} \wedge \frac{d\omega}{\omega -q}=
\frac{1}{2\pi i}\frac{\partial(\varphi,\omega)}{\partial(\psi,t)}d\psi \wedge dt.$$
On the circle $|\zeta|=1$ the values $\varphi(\zeta,t)$ coincide with the values of the meromorphic extension
$G(\zeta,t)$  and same is true for the tangential derivatives:
$$\frac{\partial \varphi}{\partial \psi}=\frac{\partial G}{\partial \psi}.$$
The functions $G$ and $\omega$ are analytic in $\Delta$ near the unit circle and smooth up to the boundary, hence
the angular derivatives express through the complex $\zeta$-derivatives:
$$ \frac{\partial G}{\partial \psi}=i\zeta  \frac {\partial G}{\partial z}, \ \frac{\partial \omega}{\partial \psi}=i\zeta  \frac {\partial \omega}{\partial z} .$$
Since 
$$d\psi=\frac{1}{i}\frac{d\zeta}{\zeta}$$
we can rewrite $\eta$ as
\begin{equation} \label{E:eta}
\eta=\frac{1}{2\pi i}(\frac{\partial_{\zeta}G}{G} \frac{\partial_t \omega}{\omega-q} -\frac{\partial_t G}{G} \frac{\partial_{\zeta} \omega}{\omega-q}) d\zeta \wedge dt. 
\end{equation}
Let $Z_j:\zeta=\zeta_1(t), \cdots, Z_N:\zeta=\zeta_N(t)$ and $P_1:\zeta=p_1(t), \cdots, P_M:\zeta=p_M(t)$ 
are  local smooth branches of zeros and poles,
in the variable $\zeta$, of the function $\varphi(\zeta,t)$ when $t$ is near $t_0 \in (0,1).$
These branches may intersect and contain boundary points $(\zeta,t) \in \partial \Delta  \times(0,1).$

From (\ref{E:eta}), the integral of $\eta$ can be computed as the double integral of the coefficient
 of the differential form $\eta$, first in 
$\zeta$ and then  in $t.$ 
Since the coefficient of $\eta$ is a meromorphic function in $\zeta,$ the first integration (in $\zeta$)
can be performed using residues . 

The singular points come from zeros and poles of the $G$ and are defined
by leading terms in the Laurent expansion.
The expansion of  $G(\zeta,t), t \in (t_0- \varepsilon,t_0 + \varepsilon), $ with respect to $\zeta$ near the zero $\zeta=\zeta_{j}(t),$
has the form 
$$G(\zeta,t)=A_j(t)(\zeta-\zeta_j(t))^{k_j}+ \mbox{higher powers} $$
and then
\begin{equation}\label{E:G_zeros}
\begin{aligned}
&\frac{\partial_{\zeta} G}{G}(\zeta,t)=\frac{k_j}{\zeta-\zeta_j(t)}+ \mbox{holomorphic terms},\\
&\frac{\partial_t G}{G} (\zeta,t) 
=\frac{A_j^{\prime}(t)(\zeta-\zeta_j(t))^{k_j}-A_j(t)k_j\zeta_j^{\prime}(t)
(\zeta-\zeta_j(t))^{k_j-1}} { A_j(t)(\zeta-\zeta(t))^{k_j}+\cdots}\\ 
&=-k_j \zeta_j^{\prime}(t) \frac{1}{\zeta-\zeta_j(t)} +\mbox{holomorphic terms}.
\end{aligned}
\end{equation}

Correspondingly, if  the Laurent series at a pole $\zeta=p_s(t)$  has the form
$$G(\zeta,t)=\frac{B_s(t)}{(\zeta-p_s(t))^{l_s}}+\cdots,$$
then
\begin{equation}\label{E:G_poles}
\begin{aligned}
&\frac{\partial_{\zeta} G}{G}(\zeta,t)=-\frac{l_s}{\zeta-p_s(t)}+ \mbox{holomorphic terms},\\
&\frac{\partial_t G}{G} (\zeta,t)=  p_s^{\prime}(t) \frac{l_s}{\zeta-p_s(t)} +\mbox{holomorphic terms}.
\end{aligned}
\end{equation}

Then for sufficiently small $\varepsilon >0$ we have from (\ref{E:eta}), (\ref{E:G_zeros}) and (\ref{E:G_poles}): 
\begin{equation}\label{E:local}
\begin{aligned}
&\int\limits_{\partial \Delta \times (t_0-\varepsilon, t_0+\varepsilon)} \eta=
\int\limits_{(t_0-\varepsilon, t_0+\varepsilon)} \sum\limits_{j=1}^N k_j 
\frac{\partial_t \omega (\zeta_j(t),t)+
\zeta^{\prime}(t)\partial_{\zeta}\omega(\zeta_j(t),t)}{\omega-q}dt \\ 
&-\int\limits_{(t_0-\varepsilon, t_0+\varepsilon)} \sum\limits_{s=1}^M p_s \frac{\partial_t \omega (p_s(t),t)+
p_s^{\prime}(t)\partial_{\zeta}\omega(p_s(t),t)}{\omega-q}dt.  
\end{aligned}
\end{equation}
In this formula, contributions from the boundary zeros come with the coefficients half of multiplicity $k_j/2.$
It is essential for the integral only if the boundary zeros fill an open curve $\zeta=\zeta(t), |\zeta(t)|=1.$

Change the variable $z=\omega(\zeta,t).$ Since 
$$dz=d\omega=\zeta_j^{\prime}(t) \partial_{\zeta}\omega d\zeta + 
\partial_t \omega dt$$  
on $Z_j$ and 
$$dz=d\omega =p_s^{\prime}(t) \partial_{\zeta}\omega d\zeta + 
\partial_t \omega dt$$ 
on $P_s,$ formula \ref{E:I(q)=}) follows from (\ref{E:local}) by the above change of variable and summing local integrals.    Lemma is proved.

Combination of Lemmas \ref{L:I(q)=0} and \ref{L:I(q)=} leads to
\begin{corollary} \label{C:corr}
For every $q \notin \overline \Omega$ holds 
\begin{equation}\label{E:Z-P=0}
\int\limits_{\omega(Z(\varphi))}\frac{dz}{z-b}-\int\limits_{\omega(P(\varphi))}\frac{dz}{z-q}=0.
\end{equation}
\end{corollary}
\begin{remark} We have used the condition for $q$ only for the integrability of $(z-q)^{-1}.$
However, this singularity is integrable and hence the same arguments show that, in fact, 
the identity (\ref{E:Z-P=0}) is true so long as the right hand side in (\ref{E:Z-P=0}) makes sense, i.e. $q \notin \omega (Z(\varphi) \cup P(\varphi)).$ 
However, for our purposes, it is enough to take  $q \notin \overline \Omega.$   
\end{remark}
{\bf Proof of Proposition \ref{P:zeros_poles}}

Consider the set $\Sigma:=\omega(Z(\varphi)) \cup \omega(P(\varphi)).$ It is the 
union of oriented curves, with the orientation  induced by the orientations on the curves of zeros and of poles of $\varphi(\zeta,t).$ The latter orientation is defined by the parameter $t \in (0,1)$. We will regard $\Sigma$ as a 1-chain, by counting multiplicities and reversing  the orientation on $\omega(P(\varphi)).$ The closed domain $\overline \Omega$ is now the closed unit disc and Corollary \ref{C:corr} reads as 
\begin{equation}\label{E:Sigma}
\int\limits_{\Sigma} \frac{dz}{z-q}=0, \ \ q \notin \overline \Delta.
\end{equation} 
But this implies that the chain $\Sigma$ is union of closed curves. 
Indeed, if $\Sigma=\Sigma_1 \cup \cdots \cup \Sigma_K$ where $\Sigma_j$ is an oriented curve with the initial point
$\alpha_j$ and the end point $\beta_j$ then 
$$\int\limits_{\Sigma}\frac{dz}{z-q}=\sum\limits_{j=1}^K (ln(\alpha_j-q)-ln(\beta_j-q))=0,$$
where a continuous  branch of $ln(z-q), z \in \Delta$ is chosen due to the condition $q \notin \overline \Delta.$ 
Since $q$ varies in the open set $\mathbb C \setminus \overline \Delta$ , the latter identity is possible only if the terms $ln(\alpha_i-q)$ and $-ln(\beta_j-q)$ pair-wise cancel each other, i.e. the initial points $\alpha_i$ and endpoints $\beta_j$ pair-wisely coincide. This means that $\Sigma $ can be represented as the union of closed curves.

By the condition, the zero set $Z(\varphi)$ contains  $N$ continuous curves, $Z_1, \cdots,Z_N$. The image chains
$\omega(Z_i) \subset \Delta, \ i=1,\cdots,N, $ start at $a$ and end at $b$. On the other hand, 
the image chains $\omega(P_j), j=1,\cdots,M$ of the continuous poles $P_1,\cdots,P_M,$ the opposite, start at $b$ and end at $a.$  The combined chain must be closed which means that the number of curves going from $a$ to $b$ and the number of curves going from $b$ to $a$ must coincide, i.e. $N=M.$ 

This completes the proof of Proposition \ref{P:zeros_poles}.

\section{Zeros and poles of the meromorphic extensions of $\overline \partial-$derivative}\label{S:dbar}

Let $\Gamma=\{\gamma_t\}, \ t \in [0,1]\}$ be a $C^1$-family  of smooth Jordan curves. 
Parameterize the corresponding domains $D_t$ by a family of conformal mappings
$$\omega(\cdot,t):\Delta \to D_t$$
smoothly depending on the parameter $t.$
Consider the Jacobian
$$Q(w,t)=\frac{\partial(\omega,\overline\omega)}{\partial(\psi,t)}, \ w=re^{\psi}.$$

If the family consists of circles, then $\omega(w,t)$ can be chosen affine in $w$:
$$\omega(w,t)=c(t)+r(t)w$$ 
and then
$$Q=-ir(c^{\prime}\overline w + \overline c^{\prime}w+2r^{\prime}w\overline w).$$
Replacing $\overline w=1/w,$ we see that $Q(w,t)$ has the following meromorphic extension from the circle $|w|=1$ into the disc $|w|<1:$
$$ \mathcal D(w,t)=-ir\frac{d(w,t)}{w},$$
where $d(w)$ is the discriminant, introduced earlier:
\begin{equation}\label{E:discrim}
d(w,t)=c^{\prime}(t)+2r^{\prime}(t)w+ \overline c^{\prime}(t)w^2.
\end{equation}
\begin{lemma}\label{L:dbar} Let $\mathcal C=\{C_t=C(c(t),r(t)), \ t \in [0,1]\}$
be a chain of circles, starting at $a$ and ending at $b \neq a.$
Let $f \in C^1(\Omega)$ and suppose that $f$ extends meromorphically
from $\mathcal C$ with the only singularities-poles at the geometric centers $c(t)$ of the circles $C_t$, of order
at most $ \nu.$ 
Let $g:=\partial_{\overline z}$ and suppose that $g$ is a nonzero function. Then
\begin{enumerate}
\item if $\nu=0$  then $g$ extends meromorphically inside any circle $C_t$ and the extension in $C_t$ has at most one simple pole in the disc $|w|<1$ and the zero at $c(t),$ of at least second order. 
\item if $\nu >0$  and  the chain $\mathcal C$ satisfies the condition (*) of Theorem \ref{T:MAIN}, then $g$ meromorphically extends inside $C_t$ with the poles at the centers $c(t),$ of 
order at most $\nu-1,$ and there is no other traveling poles for $g$. 
\end{enumerate}
\end{lemma}
\pf 
Denote $K(w,t)$ the meromorphic extension of the function $f(c(t)+r(t)w)$ in the unit disc $|w|<1.$ 
We have on the unit circle $|w|=1$: 
$$K(w,t)=f(c(a,t)+r(t)w).$$
Set  $w=e^{i\psi}$ and differentiate  this identity in $t$ and $\psi$:
\begin{equation}
\begin{aligned}
&K_t=\partial_z f (c^{\prime}+r^{\prime} w)+\partial_{\overline z}f (\overline c^{\prime}+r^{\prime}\overline w), \\
&K_{\psi}=\partial_zf (irw)+\partial_{\overline z}f (-ir \overline w)
\end{aligned}
\end{equation}
Then by the Cramer's rule
\begin{equation}\label{E:g}
g(z)=\partial_{\overline z}f(z)=
\frac{w(c^{\prime}(w,t)+r^{\prime}(w,t)w)\partial_{\psi}K(w,t) -irw^2\partial_tK(w,t)}{d(w,t)},
\end{equation}
where $z=c(t)+r(t)w.$
The identity holds on the unit circle $|w|=1$, that is when $c(t)+wr(t) \in C_t.$

Now we want to understand zeros and poles of $g$ in the disc $|w|<1$. 
First of all, the denominator $d(w,t)$ may vanish identically, but only for finite set of $t$ and to a finite order in $t,$ according to  our general conditions for the chain.
Since $g(c(t)+r(t)w)$ is continuous on $|w|=1$, the enumerator vanishes at zeros of denominator and the singularity cancels
due to the regularity condition for $f$ and hence for $g$. Thus, zeros of $d(w,t)$ in $t,$ identical in $w,$ are removable singularities in (\ref{E:g}).

Suppose  $d(w,t)$ is not the identically zero function. 
By the regularity condition, $g$ may have on $|w|=1$
only finite number of finite order zeros and so that the boundary may contribute to the zeros, but not to the poles.

Since $\partial_{\psi}K=iw\partial_{w}K$, the first term in the enumerator of  (\ref{E:g}) is divisible by $w^2.$ The second term  has the factor $w^2.$ 
If the order of the poles $\nu=0$ then $K$ is analytic in $|w|<1$ and so are the derivatives 
$K_{\psi}$ and $K_r.$  Therefore, due to the factor $w^2,$ the enumerator in (\ref{E:g}) has at $w=0$ a zero of order at most 2. Also, the quadratic polynomial $d(w,t)$ in the denominator has either both roots on the unit circle, or
has one root inside the circle $|w|=1.$ In the first case, the poles are removable because $g$ is smooth on $C_t.$
In the second case, $d(w,t)$ develops a simple pole inside the disc $|w|<1$ which corresponds to a pole inside $C_t$ for the meromorphic extension . The statement 1 is proved.

Consider the case $\nu >0.$ By the condition, $w=0$ is the pole of order at most $\nu$ for $K(w,t).$
The differentiation in parameters increase the order of the pole by 1 and hence $K_r$ and $K_{\psi}$ have poles at $w=0$ of order at most $\nu+1.$ The factor $w^2$ in the enumerator
decreases the order of the pole $w=0$ by 2, and, as the result,  the extension of $g$ in $C_t$ has a pole at the center
$c(t)$ of order at most $(\nu+1)-2=\nu-1.$  Also the denominator $d(w,t)$ may add a simple pole $p(t)$ in the disc $|w|<1.$ But the condition (*) just rules out the possibility for this pole to be traveling  
and thus the statement 2 is proved as well.

\section{End of the proof of Theorem \ref{T:MAIN}}\label{S:proof1}

Now we can finish the proof of our main result. The proof goes by induction in the order $\nu$ of poles at $c(t)$ of the
meromorphic extensions inside the circles  $C_t=C(c(t),r(t))$, under assumption that there are no other traveling poles.

Let us start with $\nu=0.$ We want to prove that $f$ is analytic in $\Omega$, i.e. 
the $\overline \partial-$ derivative $g:=\partial_{\overline z}f=0.$ Suppose that $g$ is not identically zero. Then Lemma \ref{L:dbar}, statement 1,  says that $g$ has the traveling zero at $c(t)$ of at least second order and at most one traveling simple pole. Thus, $N_g \geq 2$ while $M_g \leq 1$ and hence $N_g \neq M_g$ in contradiction with Proposition \ref{P:zeros_poles}. Therefore $g=0$ and $f$ is analytic. Notice that the condition (*) was not used.

Now, assume that $\nu>0$ and assume that the statement of Theorem \ref{T:MAIN} is proved for the order of poles $\nu-1.$
Let again $g:=\partial_{\overline z}f=0$ and suppose that $g$ is nonzero function (otherwise $f$ is analytic and there is nothing to prove). Lemma \ref{L:dbar}, statement 2, yields that the number of traveling zeros is $N_g \leq \nu-1$.
Thus, we are in the conditions of the induction assumption which yields that $g$ is polyanalytic of order $\nu-1$, i.e. 
$$\frac{\partial^{\nu}g}{\partial \overline z^{\nu}}(z)=0, \ z \in \Omega.$$
This is the same as
$$\frac{\partial^{\nu+1}f}{\partial \overline z^{\nu+1}}(z)=0,$$
that is $f$ is polyanalytic of order $\nu$ and this completes the proof.

\section{Proof of Theorem \ref{T:noregular}}
The proof of Theorem \ref{T:noregular} is similar to that of  Theorem \ref{T:MAIN} but much easier
due to absence of the discriminant set.

The initial step of the induction is the case $\nu=0$, i.e. the case of analytic extendibility into the circles $C_t.$ It was proved in \cite{T2} that if $|c^{\prime}(t)| >|r^{\prime}(t)|$ and the circles $C_t$ are not enclosed into each other then the condition of analytic extendibility implies analyticity, which is the assertion of  Theorem \ref{T:noregular} in the case $\nu=0.$  

The induction step is absolutely similar to that in the proof of Theorem \ref{T:MAIN}.
However, now the situation is even easier because the dicriminant polynomial does not produce additional poles in (\ref{E:g}) any longer. Hence by Lemma \ref{L:dbar} 
the derivative $g:=\partial_{\overline z}f$ extends in the circles $C_t$ with poles, of order at most $\nu-1,$ at the centers $c(t),$ while no other poles inside $C_t$ appear. By the induction assumption, $g$ is polyanalytic of order $\nu-1$ which is equivalent to  $f$ being polyanalytic of order $\nu$. This completes the proof.

It is worth mentioning  that  the result of \cite{T2} can not be used in the proof of Theorem \ref{T:MAIN}, as it does not cover the case of enclosed circles which is in our main focus (see Examples 1-3 in Section \ref{S:examples}).

\section {Proof of Corollary \ref{C:hyperbolic}}\label{S:proof2}

First of all, we have to prove that the condition (*) is fullfiled for the family
$\mathcal H(a,b)$ of hyperbolic circles or horicylces (depending whether the point $a$ or $b$ is inside the unit circle or belongs to it), introduced
in Examples 1-3. 

\begin{lemma}\label{L:singular}
The condition (*) of Theorem \ref{T:MAIN} holds for the family $\mathcal H(a,b)=\{H(a,r)\}\cup \{H(b,r)\}, \ 0<r<1,$
defined in Examples 1-3 in Section \ref{S:examples}.
\end{lemma}

\pf  If one of the points, say, $a$ belongs to the unit circle then the corresponding family consists of horicyles.  
In this case, we have proved in Section \ref{S:examples}. that the the discriminant set produced by the 
family $H(a,t)$ consists of one point $a$.  

Now we want to check that the same is true for the case when the point is inside the unit cricle. 
Let it be, for example,  the point $b$, i.e.  $|b| <1.$
It is convenient for simplicity  to move the points $a$ and $b$ to the real diameter by applying a suitable 
Moebius automorphism of $\Delta.$ 

We regard  $t \in (0,1) $ as the parameter for the family $H(b,t)=\{|(z-b)(1-bz)^{-1}|=t.$
The Euclidean center $c=c(t)$ and radius $r=r(t)$ of the hyperbolic circle $H(b,t)$ are
$$c(t)=b\frac{1-t^2}{1- b^2t^2}, \ r(t)=t \frac{1-b^2}{1-b^2t^2}$$
The speed $c^{\prime}(t)$ of the motion of the center and the speed $r^{\prime}(t)$ of changing the radius are:
$$c^{\prime}(t)=-2bt\frac{1-b^2}{(1-b^2t^2)^2}, \ r^{\prime}(t)= (1-b^2)\frac{1+b^2t^2}{(1-b^2t^2)^2}.$$

The discriminant equation becomes:  
$$d(w,t)=c^{\prime}(t)w^2+2r^{\prime}(t)w+c^{\prime}(t)=2\frac{1-b^2}{(1-b^2t^2)^2}(-btw^2+(1+b^2t^2)w-bt)=0.$$
It has two roots 
$$w_1=bt, \ w_2=\frac{1}{bt}.$$ 
Since $|b|<1$ and $ t \in [0,1],$ the inner root is $w_1=w_1(t)=bt.$
Then the corresponding point in the discriminant set is:
$$c(t)+r(t)w_1(t)=c(t)+r(t)bt= b\frac{1-t^2}{1-b^2t^2}+ t\frac{1-b^2}{1-b^2t^2}bt=b.$$
Thus, when $t$ runs from 0 to 1, the discriminant point remains  unmovable and is always $b$.

Thus, we arrive at the conclusion that in all cases, the discriminant set of the family $\mathcal H(a,b)$
consists only of two points:
$$\mathcal D=\{a,b\}.$$
Then  condition (*) is fullfiled for obvious reason.

{\bf End of the proof of Corollary \ref{C:hyperbolic}}.
Corollary \ref{C:hyperbolic} follows if to apply Theorem \ref{T:MAIN} to the function
$$f(z):=F(z)(1-|z|^2)^{\nu}.$$
Then $f$ is regular on $\overline \Delta$ due to the regularity $F$ in the open unit disc and the behavior at
the unit circle.

On the circle $C(c,r)$ we have:
$$1-|z|^2=1-z(\overline c + \frac{r^2}{z-c})$$
and therefore $1-|z|^2$ extends meromorphically from $C(c,r)$ with the simple pole at the geometric center $z=c.$  
On the other hand, the meromorphic extehnsion of  $1-|z|^2$ in any circle $H(c,r)$ has the simple zero at the hyperbolic center $c$. To check this, apply the conformal automorhism
$$w=\frac{z-c}{1-\overline c w}$$
of the unit disc and transform the circle $H(c,r)$ to the circle $|w|=r.$
Then the meromorphic extension of the function 
$$1-|z|^2=1-\vert \frac{w+c}{1+\overline c w}\vert ^2=\frac{(1-|c|^2)(1-|w|^2)}{|1+\overline c w|^2}$$
from the circle $|w|=r$ is
$$\frac{(1-|c|^2)(1-r^2)w}{(1+\overline c w)(w+ cr^2)}$$
and it has the simple zero at $w=0$ which corresponds to the simple zero of $1-|z|^2$ at $z=c.$ 

Thus, the factor $(1-|z|^2)^{\nu}$
cancels the poles of $f$ at $a$ and $b$ but, instead, creates poles of order $\nu$ at the geometric centers of the 
circles $H(a,r)$ and $H(b,r).$ Then the function $f$ satisfies the conditions of Theorem \ref{T:MAIN} 
for the chain $\mathcal H(a,b)$ from Example 1. 

Applying Theorem \ref{T:MAIN} , we obtain that
$$f(z)=\sum_{j=0}^{\nu}\overline z^j\widetilde{\widetilde h}_j(z),$$ 
where $\widetilde{\widetilde h}_j$ are analytic functions in  $\Delta.$
Substitution
$$\overline z=\frac{(|z|^2-1)+1}{z}$$
leads to the representation of the form
$$f(z)=\sum_{j=0}^{\nu}(1-|z|^2)^j \frac{\widetilde h_j(z)}{z^j},$$
where $\widetilde h_j$ are again analytic.
The meromorphic extension of $f$ inside any circle surrounding $0$  develops a pole
(coming from the factor $(1-|z|^2)^j$) at the geometric center and if this center in not $0$ then
the extension is analytic at $0.$  Therefore 0 can not be  a singular point and hence the functions $h_j(z):=z^{-j}\widetilde h_j(z)$ are analytic at $z=0.$ Then the  representation (\ref{E:F}) follows and this completes the proof.
\section{Concluding remarks}\label{S:concluding}
\begin{itemize}
\item
Theorem \ref{T:MAIN} is true for periodic families $\{C_t\}_{t \in S^1}$ of circles, parameterized by points on the circle $S^1$ rather than on the interval (0,1). The crucial condition, replacing the condition $a \neq b$ for shrinking chains,  is that the centers of the circles constitute a curve which is non contractible in the union of the closed discs $\overline D_t$ bounded by $C_t$ ( the condition of homological non triviality in \cite{A1,A2}).

Let us explain briefly how certain points of the proof look in the periodic case.
The condition (*)  reads now as follows: the discriminant set $S(\mathcal C)$ does not
contain a closed continuous curve that is non contractible in the union of $\overline D_t$. Then the  proof of the key statement, Proposition \ref {P:zeros_poles} is almost the same.  The proof that $I(q)=0$ in Lemma \ref{L:I(q)=} does not change because $\deg \varphi=0$ for the same reasons (the manifold $\partial \Delta \times S^1$ is closed.) In the end part of the proof of Proposition \ref{P:zeros_poles} the point $q$
is taken out of the union of the discs $\overline D_t$ but so that the curve of centers $c(t)$ has a nonzero index with respect to $q$. By the condition, such $q$ exists. Traveling zeros have the total index $N_g$ with
respect to $q$ and since $I(q)=0$ the traveling zeros must be compensated by the same number of traveling poles.
The condition (*) excludes that the discriminant set contributes in traveling poles and one concludes that $N_g=M_g.$  

\item
In the non-periodic case, the condition  
of shrinking the circles  to two points $a,\  b$  seems rather technical and presumably can be replaced by the condition that  there is no point surrounded by all the circles $C_t,$ as in \cite{A1,A2}. 
For example, one can require, as in \cite{T2}, that the initial and end circles $C_0$ and $C_1$ belong to the exterior of each other.
We preferred here to consider the case of  shrinking families not to make the proof more complicated and to cover the case of the family $\mathcal H(a,b)$ of hyperbolic circles and horicycles, described in  Section \ref{S:examples}). 
\item The equality $I(q)=0$ from Lemma \ref{L:I(q)=} which is crucial for the proof of Proposition \ref{P:zeros_poles}, in all the examples in Section \ref{S:examples} can be proved directly without using the fact that the topological degree $\deg \varphi=0.$ Each family in Examples 1-3 consists of two
foliations, by enclosed circles,  of the unit disc. The circles from the first foliation are growing,
when the parameter $t$ increases, till they reach the unit circle a the moment $t=1/2.$. 
The second foliation, parametrized by $t \in [1/2,1],$ is similar but
the circles monotonically decrease, starting from the unit circle and shrinking at the end to a point. 
Each foliation defines a curvilinear
coordinate system in the unit disc,  and these systems have the opposite orientations.
When $t$ runs from 0 to 1 then the circles cover the disc $|z| \leq 1$ twice and the substitution $z=\omega(\zeta,t)$ in (\ref{E:defI}) leads to 
$$
\begin{aligned}
&I(q):=\int\limits_{\partial \Delta \times (0,1/2)}\eta + \int\limits_{\partial \Delta \times (1/2,1)}\eta=\\
&\int\limits_{|z| \leq 1}\frac{\partial_{\overline z}g}{g} \frac{1}{z-q}d\overline z \wedge dz-\int\limits_{|z| \leq 1}\frac{\partial_{\overline z}g}{g} \frac{1}{z-q}d\overline z \wedge dz=0.
\end{aligned}
$$
\item It would be interesting to get rid of the  regularity condition in Theorem \ref{T:MAIN}.
In our proof this condition is essential, since without assumption of regularity the structure of zeros
may be complicated and the logarithmic residue (the integral (\ref{E:defI}) may be not defined at all. 
\item The Euclidean analog of Corollary \ref{C:hyperbolic} is not true. The following example is given in \cite{AG}: if $c_1,\cdots, c_l$ is an arbitrary finite set then the function $f(z)=\overline z(z-c_1) \cdots  (z-c_l)$ 
extends analytically inside any circle centered at one of the points $c_j$  but $f$ is not analytic. The crucial difference with the case of two families of concentric hyperbolic circles (Example 1 in Section \ref{S:examples}) is that in the above example the family is disconnected.

\item Moment type characterization of polyanalytic functions was given in \cite{Z1}.
There, one zero moment is required, while the testing  family is two-parametric, it consists of circles with arbitrary centers and two well chosen radii. In our case the situation is different: 
the family of circles depends on one parameter, but 
instead  vanishing of {\it all} complex moments of order $\geq \nu$ is required. 
\item Theorem \ref{T:MAIN} gives a Morera type characterization of solutions of the polyanalytic differential equation
$\frac{\partial^{\nu+1} f}{\partial \overline{z}^{\nu+1}}=0$. Various aspects of the problem of characterization of solutions of differential equations by integral conditions of moment type were  studied in \cite{Z2}, \cite{E}, \cite{AN}. 
\end{itemize} 
 
 After this article was written,  Josip Globevnik  informed the author that he has a proof
of the statement of Theorem \ref{T:MAIN} for the partial case $\nu=0$ (corresponding to the characterization of analytic functions), but, instead, under a priori assumptions of only continuity of the functions.

\section*{Acknowledgement} This work was partially supported by ISF (Israel Science Foundation), Grant 688/08.
Some of this research was done as a part of European Networking Program HCAA.

\end{document}